\input amssym.def
\input amssym
\input epsf.tex
\input epsfig.sty
\documentstyle[11pt]{article}
\begin{document}
\bibliographystyle{alpha}

\newtheorem{thm}{Theorem}
\newtheorem{defin}[thm]{Definition}
\newtheorem{lemma}[thm]{Lemma}
\newtheorem{propo}[thm]{Proposition}
\newtheorem{cor}[thm]{Corollary}
\newtheorem{conj}[thm]{Conjecture}
\newtheorem{exa}[thm]{Example}

\centerline{\LARGE \bf Quandles and Monodromy}

\vspace*{1cm}

\centerline{\parbox{2.5in}{D.\ N.\ Yetter \\ Department of Mathematics \\
Kansas State University \\ Manhattan, KS 66506}}
\vspace*{1cm}

{\small
\noindent{\bf Abstract:} We show that a variety of monodromy phenomena 
arising in geometric topology and algebraic geometry are most conveniently 
described in terms of homomorphisms from a(n augmented) knot quandle associated
with the base to a suitable (augmented) quandle associated with the fiber.  
We consider the cases of
the monodromy of a branched covering, braid monodromy and the monodromy of
a Lefschetz fibration.  The present paper is an expanded and corrected 
version of \cite{Y.Dq}}.

\section{Introduction}

Monodromy phenomena are usually modelled by group homomorphisms from the
fundamental group of a base space with a subspace of singular values
or branch points
deleted to a suitable group of automorphisms of a generic fiber.  These
homomorphism must then satisfy 
side-conditions, unnatural from the point of view of group
theory, which are imposed on the homomorphism by the geometry of the situation.
Because of these side-conditions, some authors prefer to replace the 
group homomorphism with a map from a generating set on which the 
side-conditions are more natural (cf. \cite{GoSt,Moi}).  
This latter approach, however,
has the drawback of requiring the introduction of a combinatorial equivalence
by moves as a replacement for algebraic homomorphisms.

It is the purpose of this paper to show that in the cases considered, all
of the side-conditions can be replaced by purely algebraic
conditions by describing monodromy not in terms of a group homomorphisms or
maps on generating sets, but by homomorphisms between quandles associated to
the base and fiber, or quandles augmented in the groups usually used to 
describe the monodromy.

\section{Quandles, Fundamental Quandles and Knot Quandles}

Quandles were originally introduced by Joyce \cite{Joyce1,Joyce2} 
as an algebraic
invariant of classical knots and links.  They may be regarded as an abstraction
from groups inasmuch as some of the most important examples arise by
considering a group with left and right conjugation as operations.

\begin{defin}
A {\em quandle} is a set Q equipped with two binary operations $\rhd$ and
$\unrhd$ satisfying

\[
\begin{array}{ll}
\forall x \in Q & x \rhd x = x \\
\forall x,y \in Q & (x \rhd y) \unrhd y = x = (x \unrhd y) \rhd y \\
\forall x,y,z \in Q & (x \rhd y) \rhd z = (x \rhd z) \rhd (y \rhd z)
\end{array} \]
\end{defin}

Algebraic structures satisfying the second and third axioms only have been
studied under the name ``racks'' by Fenn and Rourke \cite{FR}. Structures
satisfying the third axiom only are called ``right distributive semigroups'' by
universal algebraists.

Examples abound:  

\begin{exa}
If $G$ is any group, we can make $G$ into a quandle
by letting $x \rhd y = y^{-1}xy$ and $x \unrhd y = yxy^{-1}$.
Likewise any union of conjugacy classes in a group $G$ forms a subquandle.
\end{exa}

This example is of particular importance for the theory of quandles, as a 
representation theorem due to Joyce \cite{Joyce1} shows that all free
quandles embed into (free) groups as a disjoint union of conjugacy
classes, and thus that the universally
quantified equations holding in all quandles are precisely those holding
in all quandles of this form.

It will also be of interest in the present investigation, as the most of
the quandles considered herein can be identified with unions of conjugacy 
classes in groups arising naturally in geometric topology.

\begin{exa}
Given an $R$-linear automorphism $T$ of an $R$-module $V$, 
$V$ becomes a quandle
with $x \rhd y = T(x-y) + y$ and $x \unrhd y = T^{-1}(x-y) + y$.
\end{exa}

The following examples are of particular interest to us, as we will see a 
topological application later:

\begin{exa} \label{alternating.quandle}
Let $R$ be any commutative ring, and $X$ be a free $R$-module equipped with an 
antisymmetric bilinear form $\langle -,- \rangle:X \times X \rightarrow R$. 
(If $(R,+)$ has any two-torsion, we actually need alternating rather than just
antisymmetric.)
Then $X$ is a quandle when equipped with the operations

\[ x \rhd y = x + \langle x,y \rangle y\]

and 

\[ x \unrhd y = x - \langle x,y \rangle y \]
\end{exa}

The proof that this last example satisfies the quandle axioms is routine, but 
we indicated it to give the reader unfamiliar with quandles the flavor of such
things:

Observe that since $\langle -,- \rangle$ is alternating, we have 
$x \rhd x = x$.  Likewise by bilinearity and alternating-ness, it follows that
$\langle x+ay,y \rangle = \langle x,y \rangle$ from which the second quandle
axiom follows.  For the third, we calculate

\begin{eqnarray*}
	(x \rhd y)\rhd z & = & (x + \langle x,y \rangle y) \rhd z \\
			 & = & x + \langle x,y \rangle y + 
				\langle x + \langle x,y \rangle y, z\rangle z\\
			 & = & x + \langle x,y \rangle y + 
				\langle x,z \rangle z + \langle x,y \rangle
				\langle y,z \rangle z \\
\end{eqnarray*} 

\noindent while

\begin{eqnarray*}
	(x \rhd z)\rhd (y \rhd z) & = & 
			(x + \langle x,z \rangle z) \rhd
				(y + \langle y,z \rangle) z \\
			& = & x + \langle x,z \rangle z + 
			 \langle x + \langle x,z \rangle z, 
				y + \langle y,z \rangle z\rangle
				(y + \langle y,z \rangle z) \\
			& = & x + \langle x,z \rangle z + \\ 
			&  &\hspace*{1cm} [\langle x,y \rangle y + 
				\langle x,z \rangle \langle z,y \rangle + 
				\langle y,z \rangle \langle x,z \rangle +
				\langle x,z \rangle \langle y,z \rangle
					\langle z,z \rangle ]
				(y + \langle y,z \rangle z) \\
			& = & x + \langle x,y \rangle y + 
				\langle x,z \rangle z + \langle x,y \rangle
				\langle y,z \rangle z \\
\end{eqnarray*} 

\noindent where the equations follow from the definitions (twice), bilinearity
and alternating-ness respectively.

We will call a quandle arising in this way {\em an alternating quandle}.

We will not directly apply alternating quandles, but rather a particular
quotient quandle which exists for any alternating quandle:

\begin{exa}
Given an alternating bilinear form $\langle -,- \rangle$ on an $R$-module
$V$, the alternating quandle structure on $V$
induces a quandle structure on the space of orbits of the
action of the multiplicative group $\{ 1, -1 \}$ on $V$ by scalar 
multiplication (note: if $1 = -1$ in $R$ the action is trivial):
negating $x$ negates $x \rhd y$ and $x \unrhd y$, while negating $y$ 
leaves them unchanged (since the negation of the instance of $y$ in
the bilinear form cancels the negation of $y$ outside).
\end{exa}

We will call a quandle arising in this way {\em a reduced alternating
quandle}.\smallskip

Joyce's principal motivation in considering this structure was to provide an
algebro-topological invariant of classical knots more sensitive than the
fundamental group of the complement.

We will need the corresponding notion in arbritary dimensions.  
We consider pairs of a 
space and a subspace, equipped with a point in the complement of the subspace
$(X,S,p)$.  In particular we consider the ``noose'' or ``lollipop'':  
$(N,\{0\},2)$ where $N$ is the subspace of $\Bbb C$ consisting of union
of the unit disk and the line segment $[1,2]$ in the real axis.

By a map of pointed pairs we mean a continuous map which preserves the
base point and both the subspace and its complement.  We can then make

\begin{defin} The {\em fundamental quandle} $\Pi(X,S,p)$ of a pointed 
pair $(X,S,p)$
is the set of homotopy classes of maps of pointed pairs (where homotopies
are through maps of pointed pairs), equipped with the operations
$x \rhd y$ (resp. $x \unrhd y$) induced by appending the path
from the base point obtained by traversing $y([1,2])$, followed by
$y(S^1)$ oriented counterclockwise (resp. clockwise), followed by
traversing $y([1,2])$ in the opposite direction to the path $x([1,2])$ and
reparametrizing.
\end{defin}

For the proof that this gives a quandle structure, see \cite{Joyce1,Joyce2}.

In the case where both the space and its subspace are smooth oriented manifolds
and the subspace is of codimension 2, it is possible to identify a
particularly interesting subquandle of the fundamental quandle.

\begin{defin} The {\em knot quandle} $Q(M,K,p)$ of a pointed pair 
$(M,K,p)$, where
$M$ is a smooth manifold, $K$ a smooth embedded submanifold of codimension 2
is the subquandle of $\Pi(M,K,p)$ consisting of all maps of the noose
such that the bounding $S^1$ has linking number $1$ with $K$.  (Note: this
is in the signed sense.)
\end{defin}

Joyce \cite{Joyce1,Joyce2}
showed that the knot quandle of a classical knot determined the knot up
to orientation.

We, however, will be concerned here with 
knot quandles in general.  In particular,
in our discussion of branched coverings, we will need to consider 
knot quandles in all dimensions. In the discussions of braid mondromy
and of the mondromy of Lefschetz fibrations we will consider
the knot quandle of an oriented set of points in a surface:  
Given a (path connected) oriented surface
$\Sigma$, equipped with a finite set of points $S$, and a point $p$ not lying
in $S$, the quandle $Q(\Sigma, S, p)$ has as elements all isotopy classes
of maps of pointed pairs from the noose to $(\Sigma, S, p)$ which map the
boundary of the disk with winding number $\pm 1$ with the sign given by
the orientation of the point (always positive, except in the case of
achiral Lefschetz fibrations).

It is easy to see that there is a relationship between $Q(\Sigma, S, p)$ and
$\pi_1(\Sigma \setminus S, p)$:  an action of the fundamental group
$\pi_1(\Sigma \setminus S, p)$ on $Q(\Sigma, S, p)$ by quandle
homomorphisms is given by appending a loop representing an element of
$\pi_1$ to the initial path of the noose and rescaling.

There is, however, an more intimate relationship between $Q(\Sigma, S, p)$
and $\pi_1(\Sigma \setminus S, p)$:  

\begin{defin} \cite{Joyce1,Joyce2}
An {\em augmented quandle} is a quadruple 

\[(Q,G,\ell:Q\rightarrow G,\cdot:Q \times G \rightarrow Q) \]

\noindent where $Q$ is a quandle, $G$ is a group, $\cdot$ is a right-action of
$G$ on $Q$ by quandle homomorphisms, and the set-map $\ell$ (called
the {\em augmentation}) satisfies

\begin{eqnarray*}
	q \cdot \ell(q) & = & q\\
	\ell(q\cdot \gamma) & = & \gamma^{-1}\ell(q)\gamma 
\end{eqnarray*}
\end{defin}

\begin{propo}
For any oriented manifold $M$ with an oriented, properly embedded
codimension 2 submanifold $K$ and a point $p \in M \setminus K$, the
quadruple

\[ (Q(M,K,p), \pi_1(M\setminus K, p), \ell, \cdot) \]

\noindent where $\cdot$ is the action described above, and $\ell(q)$
is the homotopy class of the loop at $p$ which traverses the arc,
then the boundary of the disk counterclockwise, then the arc back to $p$,
is an augmented quandle.  We call the loop at $p$ just described as 
a representative for $\ell(q)$ the {\em canonical loop} of the noose
$q$.
\end{propo}

\noindent{\bf proof:} Having noted that the action of $\pi_1(M\setminus K,p)$
is by quandle homomorphisms (a fact which follows essentially by conjugation
in the fundamental groupoid---the reader may fill in the details), it remains
only to verify that the map $\ell$ satisfies the two conditions specified
in the definition of augmented quandles.
	
The first reduces to the idempotence of the quandle operation.  The second
follows from the fact that the appended loop occurs twice in the 
specification of $\ell(q\cdot \gamma)$, initially in the outgoing arc from 
$p$ with
positive orientation, and again in the incoming arc to $p$ with reversed
orientation. $\Box$ \smallskip

We also have

\begin{propo} \label{knotgen}
If $M$ is simply connected, 
the image of the augmentation, $\ell(Q(M,K,p))$ generates 
$\pi_1(M\setminus K,p)$.
\end{propo}

\noindent {\bf proof:}  This follows from van Kampen's Theorem:  killing
all of the noose boundaries kills the fundamental group, and thus the noose
boundaries generate. $\Box$ \smallskip

\section{Quandles of Cords}

In \cite{KaMa} Kamada and Matsumoto introduce a geometric construction
for a family of quandles closely related to the braid groups of surfaces.
In this section we describe their construction and adaptations of it
better suited to handling Moishezon's braid monodromy (cf. \cite{Moi}).

The constructions of this section give rise to quandles associated to
surfaces equipped with a finite set of (interior) points.  

\begin{defin}
Let $\Sigma$ be a surface, possibly with boundary, and $P \subset \Sigma$
a finite set of interior points.  A {\em cord} in $(\Sigma, P)$ is an 
embedding
of pairs (as defined in the previous section, but dropping base points) from
$c:([0,1], \{0,1\})\rightarrow (\Sigma,P)$.

Two cords are {\em equivalent} if there is an isotopy (rel boundary)
 of $\Sigma$ which fixes $P$ and carries one to the other.
\end{defin}

Equivalence classes of cords, or cords labeled with integers, will form
the elements of the quandles described in the section.  The operations
will be described using

\begin{defin}
Given a cord $\alpha$ in $(\Sigma, P)$, the {\em disk twist} around $\alpha$
is the isotopy (rel boundary) class of self-diffeomorphisms of $(\Sigma, P)$
represented by any smoothing $\Phi_\alpha$ 
of any map $\phi_\alpha$ constructed as follows:

Choose a neighborhood $N$ 
of $Im(\alpha)$ diffeomorphic to a disk, whose closure
is disjoint from $P\setminus Im(\alpha)$, and a chart identifying the
neighborhood with the disk $\{ z | \; |z| < 2\}$ and  $Im(\alpha)$ with
the interval $[-1,1]$ on the real axis.  Then $\phi_\alpha$ is given
by the identity map outside $N$, and in local coordinates by

\[
\phi_\alpha(z) = \left\{ \begin{array}{ll}
				-z & \mbox{if $|z| \leq \frac{3}{2}$}\\
				ze^{i\pi(2-\frac{2}{3}|z|)} & \mbox{if
					 $|z| > \frac{3}{2}$}\\
			\end{array} \right. \]

\end{defin}

Observe that the isotopy class of the disk twist is independent of the
choices (chart and smoothing) used in its construction, and depends only
on the equivalence class of the cord.

We can then make the following definition

\begin{defin}
The {\em quandle of cords on $(\Sigma,P)$} denoted $X(\Sigma,P)$ has
as elements the equivalence classes of cords on $(\Sigma,P)$, with operations
given on representatives by

\[ [\alpha] \rhd [\beta] = [\Phi_\beta(\alpha)] \]

\[ [\alpha] \unrhd [\beta] = [\Phi_\beta^{-1}(\alpha)] \]
\end{defin}

In the case of a disk $D$ equipped with a finite set of $n+1$ interior 
points, this quandle can be identified with the subquandle of $B_{n+1}$,
the $n+1$ strand Artin braid group under conjugation, of all conjugates
of the (positive) braid generators.

Kamada and Matsumoto \cite{KaMa}
give generators and relations for quandles of cords
in the disk (or plane) and 2-sphere.

These quandles are almost the appropriate quandles for the discussion
of Moishezon's braid monodromy \cite{Moi}. 
However, we need a slight modification:

\begin{defin}
Let $L$ be a set of non-zero integers.  An {\em $L$-cord} is a cord in
$(\Sigma, P)$ labeled with an element of $L$.  Two $L$-cords are equivalent if
they have the same label and equivalent underlying cords.

The {\em quandle of $L$-cords on $(\Sigma, P)$}, $L-C_{\Sigma,P}$ 
has as elements the
equivalence classes of $L$-cords on $(\Sigma, P)$, with operations given
on representatives by

\[ ([\alpha],l) \rhd ([\beta],\lambda) = ([\Phi_\beta^\lambda (\alpha)],l) \]

\[ ([\alpha],l) \unrhd ([\beta],\lambda) = ([\Phi_\beta^{-\lambda}(\alpha)],l)
 \]
\end{defin}

The quandles of cords $L$-cords in $(\Sigma, P)$ each have natural 
augmentations in the $|P|$-strand braid group on $\Sigma$, 
$B(\Sigma, P) = \pi_1((\Sigma^{|P|}\setminus \Delta)/{\frak S}_{|P|}, P)$.

The augmentation maps a cord $\beta$ (resp. $L$-cord $(\beta,\lambda)$) 
to the braid given by 
fixing the other
points of $P$ and moving the endpoints of the cord by the obvious isotopy
from the identity to the disk twist which at no time moves any point
through more than $\pi$ in the local polar coordinates used to describe
the disk twist (resp. a composition of $\lambda$ copies of this isotopy).  
(Notice:  all disk twists are
isotopic to the identity once one lifts the requirement that isotopies
fix $P$.)

\section{Dehn Quandles}

We now consider another geometric construction of quandles, 
related to mapping class groups of surfaces in a way weakly analogous
to the relationship between knot quandles and fundamental groups.

From \cite{Birman} we recall

\begin{defin} If $\Sigma$ is a surface, the {\em mapping class group
of} $\Sigma$ is the group $M(\Sigma) = \pi_0({\cal F}\Sigma)$, 
where ${\cal F}\Sigma$
is the group of all oriention-preserving self-diffeomorphisms of $\Sigma$,
endowed with the compact-open topology.
\end{defin}

Birman \cite{Birman} actually defines more general objects depending on
a set of distinguished points lying in $\Sigma$.  Following the usual
convention, if $\Sigma$ is of genus $g$, we denote its mapping class 
group by $M(g,0)$, the $0$ indicating the lack of distinguished points.

It is easy to verify that $M(0,0)$ is trivial. It is also well-known that
$M(1,0) \cong SL(2,{\Bbb Z})$.

Birman and Hilden \cite{BH} gave a finite presentation for $M(2,0)$.  
Building on work of McCool \cite{McC} and Hatcher and Thurston \cite{HT},
Harer \cite{Ha} gave finite presentations for the higher genus case,
which were improved by Wajnryb \cite{Wa}.

The key to approaching presentations of mapping class groups, and to
our related quandles, however, predates these developments, and is due
to Dehn \cite{Dehn}.  It depends upon a particular construction of
self-diffeomorphisms from an embedded curve:

\begin {defin} Let $\Sigma$ be an oriented surface, and $c$ a simple
closed curve lying in $\Sigma$.  $c$ then admits a bicollar neighborhood $U$.
If we identify this bicollar neighborhood with the annulus 

\[ A = \{z | 1 < |z| < 2 \} \]

\noindent in $\Bbb C$ by an orientation preserving diffeomorphism, $\phi:U\rightarrow A$ which maps $c$ to $\{z | |z|=\frac{3}{2} \}$
given in polar coordinates by $\phi = (r_\phi,\theta_\phi)$, the
self-homeomorphism $t_c^+:\Sigma \rightarrow \Sigma$ given by

\[ t_c^+(x) = \left\{ \begin{array}{ll}
			x & \mbox{if $x \in \Sigma \setminus U$} \\
			\phi^{-1}(r_\phi(x),\theta_\phi(x) + 2\pi r_\phi(x)) &
				\mbox{if $x \in U$}
		    \end{array} \right. \]

\noindent or any self-diffeomorphism
obtained by smoothing $t_c^+$ is called a {\em positive (or left-handed) Dehn
twist about $c$}. 

{\em Negative (or right-handed) Dehn twists} are defined similarly using

\[ t_c^-(x) = \left\{ \begin{array}{ll}
			x & \mbox{if $x \in \Sigma \setminus U$} \\
			\phi^{-1}(r_\phi(x),\theta_\phi(x) - 2\pi r_\phi(x)) &
				\mbox{if $x \in U$}
		    \end{array} \right. \]

\end{defin}

It is easy to see that the positive and negative Dehn twists about a curve
$c$ are inverse to each other (in the smoothed case up to isotopy).

It is well-known that the positive 
Dehn twists along isotopic simple closed curves are isotopic as 
diffeomorphisms.  Thus each isotopy class of simple closed curves determines
an element of the mapping class group.  Similarly the images of simple closed 
curves under isotopic diffeomorphisms will be isotopic.

We may thus make the following definition

\begin{defin} \label{Dehnquandle}
The {\em (chiral) Dehn quandle} $D(\Sigma)$ of an oriented 
surface $\Sigma$ is the
set of isotopy classes of simple closed curves in $\Sigma$ equipped with
the operations

\[ x \rhd y = t_y^-(x) \]

\[ x \unrhd y = t_y^+(x) \]

\noindent where by abuse of notation we use the same symbol to denote the
isotopy class and a representative curve.
\end{defin}

This quandle was originally described by Zablow \cite{Zab1,Zab2}, who did not
trouble to name it, and subsequently rediscovered by the author.

By the discussion above, it is clear that the operations described are 
independent of the choice of representing curve and define a well-defined
isotopy class of curves.
We now establish

\begin{propo}  The operations of Definition \ref{Dehnquandle} satisfy the
quandle axioms.
\end{propo}

\noindent{\bf proof:} It is clear by the discussion above 
that the second quandle axiom
is satisfied.  Likewise observe that $t_x^+$ fixes the curve $x$
up to isotopy.  Thus
the first quandle axiom is satisfied.  It thus remains only to verify the
third axiom.  This may be seen from the fact that any self-diffeomorphism of
$\Sigma$ induces an automorphism of the algebraic structure with operations
$\rhd$ and $\unrhd$, in particular $- \rhd y = t_y^+(-)$ is such an
automorphism.  $\Box$ \smallskip

In the next section we will use the Dehn quandle to encode the monodromy
of a Lefschetz fibration (whether chiral or achiral).  
To provide an alternative way of handling the case of achiral Lefschetz
fibration, we make

\begin{defin}
The {\em achiral Dehn quandle} $\tilde{D}(\Sigma)$ of an oriented
surface $\Sigma$ is the quandle with underlying set 
$D(\Sigma)\times\{+, -\}$ and quandle operation

\[ (x,\sigma) \rhd (y,-) = (x \rhd y, \sigma) \]

\[ (x,\sigma) \rhd (y,+) = (x \unrhd y, \sigma) \]
\end{defin}

We consider now Dehn quandle in 
the case of a genus one surface, where the structure of the
Dehn quandle can be completely determined.  As noted above
$M(1,0) \cong SL(2,{\Bbb Z})$.  Recall also that isotopy classes of 
essential simple closed curve are given by slopes $\frac{y}{x}$ with
$x$ and $y$ relatively prime integers (and $0$ is allowed in either place).

As noted in Casson and Bleiler \cite{CB}, elements of $SL(2, {\Bbb Z})$ 
corresponding to powers of Dehn twists are the integer matrices of 
trace $2$ and determinant $1$.  A fairly routine calculation shows that
the right-hand Dehn twist along a curve of 
slope $\frac{y}{x}$ ($\gcd(x,y) = 1$)
is given
by the matrix

\[ M_{\frac{y}{x}} = \left[ \begin{array}{cc}
			1 - xy & x^2 \\
			-y^2 & 1 + xy
	  	    \end{array} \right] \]

Observe also that this transformation from slopes to matrices is 
well-defined, being independent of the choice  of signs for $x$ and
$y$, and one-to-one:
given a matrix of the given form, $x$ and $y$ may be recovered up to sign from
the off-diagonal entries, while the diagonal entries determine the sign of
the product $xy$, and thus the sign of the slope $\frac{y}{x}$.

We can thus determine a formula for the operations of the Dehn quandle from the
observation that 

\begin{eqnarray*}
	t_{h(c)}^+ & = & h(t_c^+(h^{-1})) \;\;\; (*).
\end{eqnarray*}

Applying this fact in the case where $h$ itself is a postive Dehn twist 
gives us

\begin{exa} The Dehn quandle of the torus $D({\Bbb T}^2)$ has underlying set

\[ \{ \frac{y}{x} | x,y \in {\Bbb Z},\; \gcd(x,y) = 1 \} \cup \{I\} \]

\noindent where $I$ represents the (unique) isotopy class of contractible
simple closed curves and $\frac{y}{x}$ reprsents the isotopy class of 
essential simple closed curves of slope $\frac{y}{x}$.

The quandle operations on $D({\Bbb T}^2)$ are given by

\begin{eqnarray*}
	\frac{v}{u} \rhd \frac{y}{x} & = & 
			\frac{v - vxy + uy^2}{u + uxy - vx^2} \\
	I \rhd q & = & I \\
	q \rhd I & = & q \\
	\frac{v}{u} \unrhd \frac{y}{x} & = & 
			\frac{v + vxy - uy^2}{u - uxy + vx^2} \\
	I \unrhd q & = & I \\
	q \unrhd I & = & q 
\end{eqnarray*}

\noindent where $q$ is any element of the quandle and $x$, $y$, $u$, and
$v$ are integers with $\gcd(x,y) = \gcd(u,v) = 1$.
\end{exa}

It is easy to see that Dehn twist on contractible curves are isotopic
to the identity, and likewise that the isotopy class of contractible curves
is fixed by any Dehn twist.  The first of the
remaining two relations may be obtained by
using the equation $(*)$ above in the case where $h = t_{\frac{y}{x}}^-$ and
$c$ has slope $\frac{v}{u}$, computing the conjugate 
$M_{\frac{y}{x}}^{-1}M_{\frac{v}{u}}M_{\frac{y}{x}}$ and identifying the
numerator and denominator which give rise to the resulting matrix.  The last
remaining relation may be verified by observing that it provides the inverse
operation to that just computed.

As in the case of the fundamental and knot quandles, the Dehn quandle admits
an augmentation in the obvious related group:

\begin{propo}
There is an obvious right action of $M(\Sigma)$ on $D(\Sigma)$ by
quandle homomorphisms given by $[q]\cdot [h] = [h(q)]$.  Let 
$\ell :D(\Sigma)\rightarrow M(\Sigma)$ be given by mapping an isotopy
class of simple closed curve in $\Sigma$ it the isotopy class of the
positive Dehn twist about any of its representatives.  Then

\[ (D(\Sigma), M(\Sigma), \ell, \cdot) \]

\noindent is an augmented quandle.  We call it the {\em augmented Dehn
quandle of $\Sigma$}.

\end{propo}

\noindent {\bf proof:} The proof is routine.

As was the case with the augmented knot quandle for a simply connected
underlying manifold, so with augmented Dehn quandles we have

\begin{propo} \label{Dehngen}
The image of the augmentation $\ell(D(\Sigma))$ generates $M(\Sigma)$.
\end{propo}

\noindent{\bf proof:} This is simply a restatement of the classical theorem
that the mapping class group is generated by (right-hand) 
Dehn twists \cite{Dehn,Lickorish}. $\Box$ \smallskip

Similarly the $M(\Sigma)$ admits a right action on $\tilde{D}(\Sigma)$
by $(q,\sigma)\cdot h = (h(q),\sigma)$, and there is an augmentation
map $\tilde{\ell}:\tilde{D}(\Sigma)\rightarrow M(\Sigma)$ given by
mapping $(q,-)$ (resp. $(q,+)$ to the negative (resp. positive) Dehn
twist along $q$.  We call this the {\em augmented achiral Dehn quandle} of
$\Sigma$.

As of this writing, the structure of the Dehn quandle for higher genus 
surfaces has yet to be determined.  One thing which can be read off from the 
well-known presentation for the mapping class group for a surface $\Sigma_2$
of genus two is

\begin{propo} The Dehn quandle $D(\Sigma_2)$ of a surface of genus two admits
a quotient to a seventeen element quandle, two of whose elements act trivially 
and the other fifteen of which form the quandle of all transpositions in 
${\frak S}_6$.
\end{propo}

\noindent {\bf proof:}  First pass by the
augmentation map to the subquandle of $M(\Sigma_2)$ 
under conjugation, then to the subquandle of ${\Bbb Z}/10 \times {\frak S}_6$
under the quandle map induced by the group homomorphism which maps the 
generator $\zeta_i$ to $(1,(i\; i+1))$. The image is then the subset
$ \{ (0,e), (2,e) \} \cup \{ (1,(a\; b)) | 1 \leq a < b \leq 6 \}$. The element $(2,e)$ is the image of (any of) the product(s) of twelve Dehn twists about non-separating curves which give a Dehn twist about a separating curve, all of
which become trivial in the quotient to ${\frak S}_6$ and map to 2 in the
quotient to ${\Bbb Z}/10$.  This set
is readily verified to be closed under conjugation, which induces the 
quandle structure describe in the proposition.
 $\Box$
\smallskip

It is also possible in general to find interesting quotients of 
chiral and achiral Dehn quandles 
by considering the (reduced)
alternating quandle associated to $H_1(\Sigma,R)$ with the 
intersection form, where $R$ is any quotient of $\Bbb Z$.  We call the
alternating
quandle associated to the intersection form
the {\em $R$-homology quandle of} $\Sigma$ and denote it by
$HQ_R(\Sigma)$, omitting the $R$ when $R = {\Bbb Z}$. 
(As an aside, by the same construction, we can put a
quandle structure on $H_{2n+1}(X,R)$ for $X$ any $4n+2$ manifold.)

Since the Dehn quandle $D(\Sigma)$
has as elements isotopy classes of {\em unoriented} 
simple closed curves, they can be more naturally related to the
reduced alternating quandle associated to the interection form, which
we call the {\em $R$-homology Dehn quandle of} $\Sigma$ and denote
by $HD_R(\Sigma)$, as before
omitting the subscript $R$ when $R = {\Bbb Z}$. 

Any unoriented simple closed 
curve represents an element of $HD_R(\Sigma)$, with isotopic simple closed
curves representing the same element.  We thus have a map 
$D(\Sigma)\rightarrow HD_R(\Sigma)$ for any surface $\Sigma$.

To see that this map is a quandle map we must
relate the geometric construction of the operations in $D(\Sigma)$ to the
algebraic construction of the operation on $HD_R(\Sigma)$ from the 
intersection form.  Consider a pair $a,b$ of unoriented simple closed 
curves in an oriented
surface $\Sigma$.  Depending
on how they are oriented, their intersection number (if it is non-zero) may be
given either sign.  Since we are really 
concerned with isotopy classes of curves, 
we may assume the curves intersect tranversely.  

Now choose an orientation on $a$.  We may induce an orientation on $b$ as 
follows:  orient $b$ so that at each intersection point, the ``turn right'' 
rule defining a right-handed Dehn twist about $b$ causes the curve 
representing $a \rhd b$ in $D(\Sigma)$
to traverse $b$ with the same sign as the intersection point.

The curve representing $a \rhd b$ in $D(\Sigma)$, oriented to agree
with the orientation on $a$, then represents the 
homology class $a + \langle a,b \rangle b$ in $H_1(\Sigma,R)$.
Passing to the quotient $HD_R(\Sigma)$ then removes any dependence
on orientation, and we see that the map carrying a simple closed
curve to the $\{\pm 1\}$-orbit of its homology class is a quandle
homomorphism.

One thing which should be observed it that for genus 1, $D({\Bbb T}^2) \cong
HD({\Bbb T}^2)$ since each homology class is represented by a unique isotopy 
class of oriented simple closed curve.  In higher genus, $HD(\Sigma)$ will
be a proper quotient of $D(\Sigma)$, as different isotopy classes of curves
can represent the same homology class.  For example, both a curve which bounds
a disk and a curve which separates a surface of genus two 
into two genus one surfaces with
boundary are both null-homologous, but they represent different
isotopy classes.

It might naively be thought that the achiral Dehn quandle should have an 
an analogous relationship with the homology quandle with the signs 
in the pairs defining the elements becoming orientations on the curves.  This,
however, is not the case:  $(b,-)$ and $(b,+)$ act differently on elements of 
the achiral Dehn quandle, while in the homology quandle $a \rhd b = a \rhd -b$.

\section{Monodromy}

\subsection{Branched Coverings}

The simplest and most classical example of monodromy phenomena we will consider
is that of the monodromy of a branched covering space.  The branch set $S$ is a
codimension two subspace.  We will consider the case in which it is 
nonsingular.  Classically it
is described by considering the homomorphism from the knot group of the 
singular set, $\pi_1(B\setminus S, p)$ to the group of permuations of 
the (generic) fiber over $p$, ${\frak S}_d$, where $d$ is the number of sheets.
(cf. for example \cite{Fox,IoPi}).

However, as the local model of a branch point is given by the self-map of
$D^2 \times D^{n-2}$ by $(z,\vec{x}) \mapsto (z^k,\vec{x})$, not all group
homomorphisms arise:  only those sending meridians of $S$, 
considered as elements of $\pi_1(B\setminus S, p)$, to cyclic permutations
actually arise.  It is also common to consider {\em simple} branched coverings 
in which the monodromy of a meridian is restricted to be a transposition.

Both conditions are easily imposed by considering the monodromy as a quandle
homomorphism rather than a group homomorphism:

\begin{defin}
The {\em quandle monodromy} of a $d$-sheeted
branched covering with base $B$ and branch locus $S \subset B$ is the
quandle homomorphism $\mu: Q(B,S,p) \rightarrow {\frak C}_d$, given by mapping
any noose to the monodromy around its boundary, where ${\frak C}_d$ is the subquandle
of ${\frak S}_d$ under conjugation consisting of non-trivial
cyclic permuations.  For simple branched covers, the monodromy may be 
considered as taking values in ${\frak T}_d$, the subquandle of ${\frak S}_d$
under conjugation consisting of transpositions.
\end{defin}

In the case where $B$ is simply connected, this suffices to recove the 
monodromy in the classical sense, since the meridians generate the 
knot group of the singular set, and thus the branched covering.

For a more general base, it is still necessary to consider the classical
notion of monodromy, but here the restriction on the monodromy of
meridians can be imposed algebraically rather than combinatorially by
considering the augmentations from $Q(B,S,p)$ to $\pi_1(B\setminus S, p)$
and from ${\frak C}_d$ or ${\frak T}_d$ to ${\frak S}_d$.

\begin{defin}
The {\em augmented quandle monodromy} of a branched covering is the
homomorphism of augmented quandles whose components are the 
quandle monodromy and the monodromy of the branched covering.
\end{defin}

\subsection{Braid Monodromy}

Our second example, Moisheson's braid monodromy, takes a bit more 
description.  
It is an invariant of complex curves in ${\Bbb C}{\Bbb P}^2$, although
it can be adapted to more general surfaces in ${\Bbb C}{\Bbb P}^2$, and
when used as an invariant of the branch locus, together with the type
of monodromy just discussed, is important in the theory of symplectic
$4$-manifolds (cf. \cite{AuKa}).

Given a complex projective plane curve $V$, 
one can change coordinates so that the
curve does not pass through $[0:0:1]$. The curve then lies in the tautological
line bundle over the exceptional locus 
${\Bbb C}{\Bbb P}^1 = \{ [x:y:0] | (x,y) \in {\Bbb C}\setminus \{(0,0)\} \}$.

The inverse image of a generic point under the projection is then a set of
$d$ points, where $d$ is the degree of the curve.  At a finite set of points,
however, the curve is either tangent to the fiber or singular, and the 
inverse image will have fewer points.

Now, if one travels in a loop in the base ${\Bbb C}{\Bbb P}^1$ around a
singular point (one where the curve is either singular or tangent to the 
fiber), the monodromy is of the curve as an embedded object is an element
of the $d$-th Artin braid group, $B_d$.  

Consideration of local models show that the monodromy around a singular
point is restricted to certain conjugacy classes of $B_d$ by the geometry
of the singularity:  points of tangency have conjugate of the braid generators
as monodromy, nodes have conjugates of squares of generators, cusps have 
conjugates of cubes of generators, and so on.

Selecting a non-singular fiber, and choosing coordinates 
on ${\Bbb C}{\Bbb P}^1$
so that it is the point at infinity, one then has an affine plane curve 
fibered over ${\Bbb C}$.

We will restrict our attention to curves with at worst cuspidal singularities,
and assume (by a small perterbation if necessary) that the tangencies and
singularities all occur in different fibers.

Moishezon then makes

\begin{defin}The {\em braid monodromy} of the (affine)
plane curve $V$ is the group
homomorphism 

\[ \mu:\pi_1({\Bbb C}\setminus S,p)\rightarrow B_d \] 

\noindent obtained by
identifying the braid group of the fiber over a non-singular point $p$ with
$B_d$ and mapping generating loops to the induced monodromy.  Here $S$ is the
set of singular values for the projection.
\end{defin}

There are two restrictions on the homomorphisms thus arising.  The first, 
arising from the geometry near the singular fibers, is readily handled by
considering quandles.

\begin{defin}The {\em quandle monodromy} of an affine 
cuspidal plane curve $V$ is the 
quandle homomorphism from the knot quandle $Q({\Bbb C}, S, p)$ to 
$\{1,2,3\}-C_(\pi^{-1}(p), \pi|_V^{-1}(p))$, the
quandle of $\{1,2,3\}$-cords on $(\pi^{-1}(p), \pi|_V^{-1}(p))$
obtained by mapping each noose to the $\{1,2,3\}$-cord whose action
describes the monodromy around the noose boundary.

As has been observed, each of the quandles admits an augmentation:  from
the knot quandle to the fundamental group of ${\Bbb C}\setminus S$ and from
the quandle of $L$-cords to the braid group of the plane $\pi^{-1}(p)$.
The {\em augmented quandle monodromy} of a cuspidal affine plane curve is
the augmented quandle homomorphism with the quandle monodromy and braid
monodromy as components.
\end{defin}

The second restriction required in the projective case, however, 
lives more comfortably at the level of
groups:  because the tautological bundle is non-trivial, there is a monodromy
around the (non-singular) point at infinity.  In particular, travelling
around the point at infinity, the sheets undergo a full twist (corresponding
to the twist in the line bundle).  Thus certain suitable products of the
monodromies around the singular fibers must be the full-twist in $B_d$, usually
denoted $\Delta^2$.

The inclusion $\iota:{\Bbb C} \rightarrow {\Bbb C}{\Bbb P}^1$ 
induces quotient maps
$q:Q({\Bbb C},S,p)\rightarrow Q({\Bbb C}{\Bbb P}^1,S,p)$, and quotient maps
on the fundamental groups $\pi_1(\iota):\pi_1({\Bbb C}\setminus P)\rightarrow
\pi_1({\Bbb C}{\Bbb P}^1\setminus P)$.  The kernel of this group homomorphism
is free on a single generator represented by a counterclockwise loop $\lambda$
in ${\Bbb C}$ which has $S$ lying in the bounded region.
Because of the behavior of the monodromy at infinity, we consider also an
augmentation into the 
quotient of the braid group, 
$\delta:B_{|\pi^{-1}(p)|}\rightarrow B_{|\pi^{-1}(p)|}/<\Delta^2>$.

From this, given a projective cuspidal plane curve, lying in the tautological
bundle over the exceptional locus, we obtain a commutative square of 
augmented quandle homomorphisms 

\begin{figure}
\setlength{\unitlength}{3947sp}%
\begingroup\makeatletter\ifx\SetFigFont\undefined
\def\x#1#2#3#4#5#6#7\relax{\def\x{#1#2#3#4#5#6}}%
\expandafter\x\fmtname xxxxxx\relax \def\y{splain}%
\ifx\x\y   
\gdef\SetFigFont#1#2#3{%
  \ifnum #1<17\tiny\else \ifnum #1<20\small\else
  \ifnum #1<24\normalsize\else \ifnum #1<29\large\else
  \ifnum #1<34\Large\else \ifnum #1<41\LARGE\else
     \huge\fi\fi\fi\fi\fi\fi
  \csname #3\endcsname}%
\else
\gdef\SetFigFont#1#2#3{\begingroup
  \count@#1\relax \ifnum 25<\count@\count@25\fi
  \def\x{\endgroup\@setsize\SetFigFont{#2pt}}%
  \expandafter\x
    \csname \romannumeral\the\count@ pt\expandafter\endcsname
    \csname @\romannumeral\the\count@ pt\endcsname
  \csname #3\endcsname}%
\fi
\fi\endgroup
\begin{picture}(6312,2514)(601,-2290)
\thinlines
\put(2026,-961){\vector( 1, 0){2100}}
\put(1576,-2161){\vector( 1, 0){2475}}
\put(3676,-361){\vector( 1, 0){1875}}
\put(3301,-1561){\vector( 1, 0){2025}}
\put(4726,-2011){\vector( 0, 1){825}}
\put(1051,-1936){\vector( 0, 1){750}}
\put(5776,-1336){\vector( 0, 1){750}}
\put(2776,-1411){\vector( 0, 1){825}}
\put(2326,-511){\vector(-3,-1){675}}
\put(2401,-1711){\vector(-3,-1){675}}
\put(5626,-511){\vector(-2,-1){600}}
\put(5551,-1711){\vector(-3,-2){450}}
\put(676,-1036){\makebox(0,0)[lb]{\smash{\SetFigFont{12}{14.4}
	{rm}$L-C_{{\Bbb C},\pi^{-1}(p)}$}}}
\put(2476,-436){\makebox(0,0)[lb]{\smash{\SetFigFont{12}{14.4}
	{rm}$L-C_{{\Bbb C},\pi^{-1}(p)}$}}}
\put(601,-2236){\makebox(0,0)[lb]{\smash{\SetFigFont{12}{14.4}
	{rm}$Q({\Bbb C}{\Bbb P}^1,S,p)$}}}
\put(2476,-1636){\makebox(0,0)[lb]{\smash{\SetFigFont{12}{14.4}
	{rm}$Q({\Bbb C},S,p)$}}}
\put(4276,-2236){\makebox(0,0)[lb]{\smash{\SetFigFont{12}{14.4}
	{rm}$\pi_1({\Bbb C}{\Bbb P}^1\setminus S,p)$}}}
\put(5476,-1636){\makebox(0,0)[lb]{\smash{\SetFigFont{12}{14.4}
	{rm}$\pi_1({\Bbb C}\setminus S,p)$}}}
\put(4276,-1036){\makebox(0,0)[lb]{\smash{\SetFigFont{12}{14.4}
	{rm}$B_{|\pi^{-1}(p)|}/<\Delta^2>$}}}
\put(5701,-436){\makebox(0,0)[lb]{\smash{\SetFigFont{12}{14.4}
	{rm}$B_{|\pi^{-1}(p)|}$}}}
\put(6901,-886){\vector( 0, 1){825}}
\put(6676,-1186){\vector(-2,-1){450}}
\put(6676, 14){\vector(-2,-1){450}}
\put(6751, 89){\makebox(0,0)[lb]{\smash{\SetFigFont{12}{14.4}{rm}$<\Delta^2>$}}}
\put(6751,-1111){\makebox(0,0)[lb]{\smash{\SetFigFont{12}{14.4}{rm}$<\lambda>$}}}
\end{picture}

\caption{The square of augmented quandle homomorphisms giving the braid 
monodromy of a projective plane curve with the induced map on group kernels}
\end{figure}
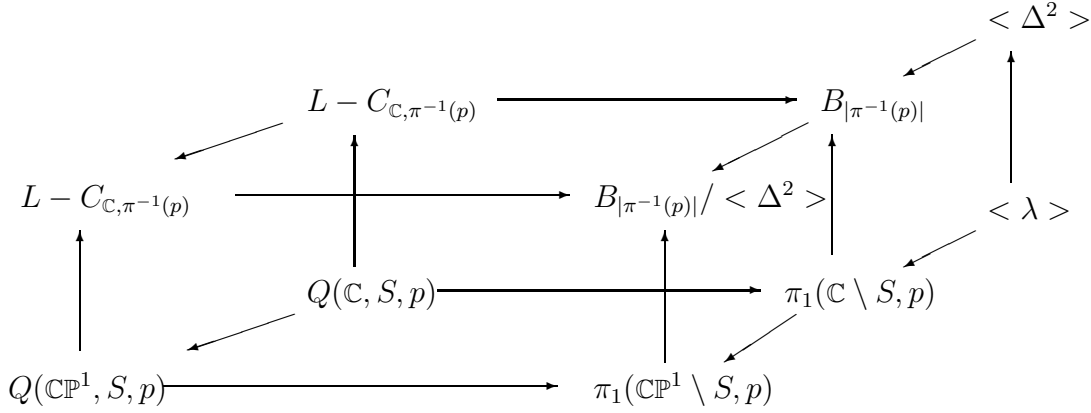 

The restriction on the monodromy at infinity then becomes the requirement
that the induced map on the kernels of the group homomorphisms maps $\lambda$
to $\Delta^2$.  Dropping this restriction will give monodromies corresponding 
to surfaces in various line bundles over ${\Bbb C}{\Bbb P}^1$---in particular,
those squares in which the induced map on kernels maps $\lambda$ to 
$(\Delta^2)^k$ will correspond to cuspidal surfaces in the line bundle with
first Chern class $k$.

We may then make

\begin{defin}The {\em augmented quandle monodromy} of a projective plane curve
is the augmented quandle monodromy of the associated affine plane 
curve.  It will necessarily satisfy the condition that the element $\lambda$
of $\pi_1({\Bbb C}\setminus S)$ is mapped to $\Delta^2$.
\end{defin}

\subsection{Lefschetz Fibrations}

Our third example is the monodromy of a Lefschetz fibration.

We briefly recall the relevant facts about Lefschetz fibrations, following
Gompf and Stipsicz \cite{GoSt}:

\begin{defin}
A {\em Lefschetz fibration} of a smooth, compact oriented 4-manifold $X$ 
(possibly with boundary) is a
smooth map $f:X\rightarrow \Sigma$, where $\Sigma$ is a compact connected 
oriented surface, $f^{-1}(\partial \Sigma) = \partial X$ and such that
each critical point of $f$ lies in the interior of $X$ and has an oriented 
local coordinate chart modelled (in complex coordinates) by 
$f(z,w) = z^2 + w^2$.

We moreover require that each singular fiber have a unique singular point.

An {\em achiral Lefschetz fibration} is defined similarly, except that the
prescribed
local coordinate chart at singularities may be orientation reversing.
\end{defin}

Now the generic fiber $F$ of $f$ is a compact, canonically oriented surface.
The genus of $F$ is called the genus of the fibration $f$.

As is pointed out in \cite{GoSt}, the choice of a regular point of the
fibration $p \in \Sigma$ and an identification of the fiber over $p$ with
a standard surface $F$ of the appropriate genus gives rise to a group
homomorphism $\Psi:\pi_1(\Sigma \setminus S, p)\rightarrow M(F)$, where $S$
is the set of critical values of $f$, called the
{\em monodromy representation} of $f$.   

In the case of genus $g \geq 2$, this group homomorphism completely determines
the structure of the Lefschetz fibration by a theorem of Matsumoto \cite{Mats}.
There are, however, restrictions on which group homomorphisms can occur.  
In particular the image of any loop linking exactly one critical value
with linking number one must be a positive Dehn twist about the
vanishing cycle of the singularity---the simple closed curve which collapses
to a point at the singular point \cite{GoSt}.

Due to the awkwardness of imposing such a condition while trying to work
in a group theoretic context, when discussing Lefschetz fibrations over the
disk $D^2$ and the sphere $S^2$, Gompf and Stipsicz \cite{GoSt} work instead
with {\em the monodromy} of the fibration:  the $|S|$-tuple of Dehn twists
given by a family of generating loops each of which links a single
critical value with linking number one.  

This, then, has the drawback that the $|S|$-tuple is determined only up to
an overall 
conjugation by an element of $M(F)$, cyclic permutation, and combinatorial
moves given by swapping two of the Dehn twists while conjugating one
of them by its partner in a suitable sense.

Both drawbacks---the use of geometric
side-conditions in what would otherwise be the the purely
group theoretic setting monodromy representations,
and the ambiguity of definition inherent in the notion of the
mondromy, are removed by considering

\begin{defin}
The {\em quandle monodromy} of a Lefschetz fibration $f:X\rightarrow \Sigma$
with critical set $S \subset \Sigma$,
relative to a regular point $p$, is the quandle homomorphism 

\[ \mu:Q(\Sigma, S, p) \rightarrow D(F) \]

\noindent given by mapping each element of $Q(\Sigma, S, p)$ to the 
monodromy around the canonical loop of any representing noose.  Here all 
points of $S$ are given the positive orientation.

The {\em quandle monodromy} of an achiral Lefschetz fibration is defined
in the same way, except that the points of $S$ are oriented positively if the
local chart modelling the corresponding singularity is orientation preserving,
and negatively if it is orientation reversing.

The {\em augmented quandle monodromy} of a Lefschetz fibration (resp. achiral
Lefschetz fibration) 
$f:X\rightarrow \Sigma$ relative to $p$ is the map of augmented quandles
$(\mu, \Psi)$, where $\mu$ is the quandle monodromy and $\Psi$ is the
monodromy representation.
\end{defin}

We then have

\begin{propo}
The quandle monodromy of a Lefschetz fibration determines the monodromy 
representation.  Conversely the monodromy 
representation
determines the quandle monodromy.
\end{propo}

\noindent {\bf proof:} The fundamental group of the complement of the
singular set in the base is generated by the image of the knot quandle under
the augmenation, while the mapping class group of the fiber is generated by the
image of the Dehn quandle under the augmentation.  Thus the quandle monodromy
induces the monodromy representation.  Conversely, the monodromy representation
satisfies the side condition that positively oriented
noose boundaries (or noose boundaries oriented according to the sign
of the singular point in the achiral case) are mapped to positive
Dehn twists, and thus the restriction of the monodromy representation to
appropriately oriented noose boundaries is the quandle monodromy.
$\Box$

This then yields the following:

\begin{thm}
The isomorphism type of the augmented quandle monodromy determines the
isomorphism class of any Lefschetz fibration of genus $g \geq 2$.  Moreover,
if $g \geq 2$ and the base $\Sigma$ is $D^2$ or $S^2$, the isomorphism 
class of the quandle
monodromy determines the isomorphism class of the Lefschetz fibration.
\end{thm}

\noindent{\bf proof:} The first statement follows {\em a fortiori}
from the theorem of 
Matsumoto \cite{Mats}.   The second statement follows from the first,
Propositions
\ref{knotgen} and \ref{Dehngen}, and the fact that in either case
$\pi_1(\Sigma \setminus S, p)$ is free. $\Box$ \smallskip

Observe that the first statement of this formulation includes
the restriction on which homomorphisms 
$\Psi:\pi_1(\Sigma \setminus S, p)\rightarrow M(F)$ actually occur as
an algebraic rather than combinatorial condition.

In the case of $S^2$, the second statement has an analogous deficiency
to the classical formulations.
Not all quandle homomorphisms extend to
augmented quandle homomorphisms: a suitably ordered product of the
Dehn twists (images of curves under the augmentation) must be the identity
in $M(F)$.

In the case of Lefschetz fibrations over the disk $D^2$ or the sphere $S^2$, 
we have

\begin{propo}
The quandle monodromy determines the monodromy up to equivalence.  Conversely,
the monodromy determines the quandle monodromy.
\end{propo}

\noindent{\bf proof:} Given the quandle monodromy, the image of any 
minimal generating
set of nooses for the knot quandle of the base, when ordered by some linear 
restriction of the cyclic ordering induced on the generators by their crossing
the boundary of a sufficiently small disk neighborhood of the base point
is a monodromy in the
sense of \cite{GoSt}. $\Box$ \smallskip

We can also define another type of quandle mondromy in the achiral
case:

\begin{defin}
The {\em achiral quandle monodromy} of an achiral Lefschetz fibration
is the quandle homomorphism from $Q(\Sigma,S,p)$,
where all points of $S$ are oriented positively,  to $\tilde{D}(\Sigma)$,
which assigns
to each element of $Q(\Sigma,S,p)$ the monodromy around the boundary
of its noose.
\end{defin}

\section{Prospects for Quandle Invariants of Monodromy}

The foregoing suggests that a fruitful approach to studying various monodromy
phenomena could be begun by finding invariants of
quandle homomorphism and augmented quandle homomorphisms which are effectively
computable from a presentation by generators and relations.

We briefly outline several places where one might begin:

\begin{itemize}
\item Simple counting invariants:  count the number of homomorphisms of
(augmented) quandle maps (that is,
commuting squares of (augmented) quandle maps) 
from the (augmented) quandle monodromy
to a fixed (augmented) quandle map between finite (augmented) quandles.
Variants of this include counting factorizations of a fixed quandle map
from $Q(\Sigma, S, p)$ to a finite quandle through the quandle monodromy.
\item Quandles map valued invariants:  Joyce \cite{Joyce1} 
considered quandles
satifying additional axioms (e.g. involutory quandles where $\rhd = \unrhd$,
and abelian quandles which satisfy $(w \rhd x) \rhd (y \rhd z) =
(w \rhd y) \rhd (x \rhd z)$).  We may consider the induced map between
(universal) quotient quandles as an invariant of the branched covering,
plane curve or Lefschetz fibration.

Similarly, in the last case,
the map $\eta:Q(\Sigma, S, p) \rightarrow HD(F,R)$, the 
``$R$-homology
quandle monodromy'' is plainly an invariant of the Lefschetz fibration.
This particular invariant, being constructed out of homology and
intersection theory, seems likely to have some geometric significance.

\item Invariants based on the quandle (co)homology of Carter, Jelsovsky,
Kamada, Langford and
Saito \cite{CJKLS,CJKS}:  this structure may be considered in 
two ways---first as
a variant of counting invariants in their guise as counting ``colorings'', and
second homologically: the quandle monodromy giving rise to a (co)chain map
between the quandle (co)chain complexes, the (co)homology of whose cone is
then an invariant of the branched covering, projective curve
or Lefschetz fibration.

Geometric interpretation of this latter invariant would then depend upon
understanding the geometric significance of the quandle (co)homology of 
the variously geometrically described quandles.
\end{itemize}

Pursuit of any of these is beyond the scope of the present work.

\bigskip
\noindent{\Large \bf Acknowledgements}
\medskip

The author wishes to thank David Auckley, Pavel Etingof and Yukio Matsumoto 
for corrections and helpful suggestion.  The author is supported by
NSF Grant \# DMS-9971510.
\bigskip

\begin{flushleft}
\bibliography{Book}
\end{flushleft}

\end{document}